\documentclass[final]{siamart}

\usepackage{amsmath}
\usepackage{amssymb}
\usepackage{blkarray}
\usepackage{bm}
\usepackage{color}
\usepackage{tabularx}
\usepackage{caption}
\usepackage{multirow}
\usepackage[margin=1in]{geometry}
\usepackage[algo2e,ruled,vlined]{algorithm2e}
\usepackage{subcaption}
\usepackage{mathtools}
\usepackage{url}

\newcommand{\norm}[1]{\left\lVert#1\right\rVert}
\usepackage{amsfonts}
\usepackage{graphicx}
\usepackage{epstopdf}
\usepackage{algorithmic}
\ifpdf%
  \DeclareGraphicsExtensions{.eps,.pdf,.png,.jpg}
\else
  \DeclareGraphicsExtensions{.eps}
\fi
\usepackage{amsopn}

\usepackage{booktabs}

\newcommand{\Sat}[1]{\ensuremath{S_{#1}}}
\newcommand{\Pres}[1]{\ensuremath{P_{#1}}}

\newcommand{\Density}[1]{\ensuremath{\rho_{#1}}}
\newcommand{\DensityComp}[2]{\ensuremath{\rho_{#1}^{#2}}}

\newcommand{\Viscosity}[1]{\ensuremath{\mu_{#1}}}
\newcommand{\Mobility}[1]{\ensuremath{\lambda_{#1}}}
\newcommand{\Porosity}{\ensuremath{\phi}}
\newcommand{\DarcyV}[1]{\ensuremath{\bm{q}_{#1}}}
\newcommand{\AbsPerm}{\ensuremath{\bm{K}}}
\newcommand{\Gravity}{\ensuremath{\bm{g}}}

\newcommand{\RelPerm}[1]{\ensuremath{k_{r#1}}}

\newcommand{\Flux}[1]{\ensuremath{\psi ^{#1}}}
\newcommand{\MolecularDiff}[2]{\ensuremath{D_{#1}^{#2}}}
\newcommand{\DiffusionFlux}[2]{\ensuremath{j_{#1}^{#2}}}

\newcommand{\MolarMass}[1]{\ensuremath{M^{#1}}}
\newcommand{\HenryConstant}{\ensuremath{H}}
\newcommand{\IdealGasConstant}{\ensuremath{R}}
\newcommand{\Temp}{\ensuremath{T}}

\newcommand{\TheTitle}{%
  Semi-smooth Newton methods for nonlinear complementarity formulation of compositional two-phase flow in porous media
}

\newcommand{\TheShortTitle}{%
  Semi-smooth Newton methods for multiphase flow in porous media
}

\newcommand{\TheAuthors}{Quan M. Bui, Howard Elman}

\title{{\TheTitle}}%\thanks{This work was funded by the Fog Research Institute under contract no.~FRI-454.}}

\author{
  Quan M. Bui\thanks{Corresponding author. Applied Math, Stats, and Scientific Computation, University of Maryland, College Park, MD (\email{mquanbui@math.umd.edu}).}
  \and
  Howard C. Elman\thanks{Department of Computer Science and Institute for Advanced Computer Studies, University of Maryland, College Park, MD (\email{elman@cs.umd.edu}).}
}

\newcommand{\TheFunding}{%
  This work is funded by the U. S. Department of Energy Office of Advanced Scientific Computing Research, Applied Mathematics program, under Award Number DE-SC0009301\@.
}

\title{{\TheTitle}\thanks{\TheFunding}}
\headers{\TheShortTitle}{\TheAuthors}
\ifpdf%
\hypersetup{%
  pdftitle={\TheTitle},
  pdfauthor={\TheAuthors}
}
\fi

\begin{document}

\maketitle

% ---------------------------------------------
% ---------------------------------------------

\begin{abstract}
Simulating compositional multiphase flow in porous media is a challenging task, especially when phase transition is taken into account. The main problem with phase transition stems from the inconsistency of the primary variables such as phase pressure and phase saturation, i.e. they become ill-defined when a phase appears or disappears. Recently, a new approach for handling phase transition has been developed, whereby the system is formulated as a nonlinear complementarity problem (NCP). Unlike the widely used primary variable switching (PVS) method which requires a drastic reduction of the time step size when a phase appears or disappears, this approach is more robust and allows for larger time steps. One way to solve an NCP system is to reformulate the inequality constraints as a non-smooth equation using a complementary function (C-function). Because of the non-smoothness of the constraint equations, a semi-smooth Newton method needs to be developed. In this work, we 
%evaluate the performance of a semi-smooth Newton method for two C-functions: the minimum and the Fischer-Burmeister functions. We also propose a new \textcolor{blue}{inexact} Newton method which employs a smooth version of the Fischer-Burmeister function as the C-function. 
consider two methods for solving NCP systems used to model multiphase flow: (1) a semi-smooth Newton method for two C-functions: the minimum and the Fischer-Burmeister functions, and (2) a new inexact Newton method based on the Jacobian smoothing method for a smooth version of the Fischer-Burmeister function.
We show that the new method is robust and efficient for standard benchmark problems as well as for realistic examples with highly heterogeneous media such as the SPE10 benchmark.
\end{abstract}

\begin{keywords}
  % 7. Keywords that describe the paper
  semi-smooth Newton, two-phase flow, nonlinear complementarity problem, phase transition, porous media.
\end{keywords}

\section{Introduction}\label{sec:intro}
Multiphase flow is a critical process in a wide range of hydrodynamic phenomena, including carbon sequestration, reservoir simulation, and groundwater remediation. For simulation, it would be ideal to have a complete knowledge of the state and composition of the fluid phases in the flow. However, this is not an easy task given the complex physics involved. Some of the most important effects that need to be taken into consideration include capillarity, miscibility, and especially phase transitions. If not handled correctly, these effects can introduce nonphysical numerical oscillations in computational solutions of the strongly nonlinear system of partial differential equations. 

Phase transitions have posed a major challenge for multiphase, multi-component models since the 1980s. If not handled correctly, they can cause numerical oscillations in solutions of these models, making such solutions physically inconsistent and unusable.  There have been many attempts to address the problems with phase transitions and determine the correct local thermodynamic state for compositional multiphase flow. In general, most of these can be classified into two common classes of methods: flash calculation \cite{Acs85,Chen97,Chen06,Coats80,Michelsen82,Young83} and primary variables switching (PVS) \cite{Forsyth95,WuForsyth01}. Flash calculation computes the local thermodynamic state from the overall mass of the individual components. While this method is stable with regard to determining the thermodynamic state, it tends to be inefficient because it requires solution of a large nonlinear system of equations at each time step (in addition to solution of the linearized systems) to recover all the thermodynamic quantities of interest. The second class, PVS, involves adapting the primary variables to the thermodynamic constraints locally. The idea is that whenever phase transitions occur, physical variables that are physically inconsistent (indicated by negative saturation, for example) are switched to well-defined quantities. The governing equations related to those variables are also modified accordingly. Although this approach is locally more efficient than flash calculations, it suffers from irregular convergence behavior in the nonlinear solve, which is typically addressed by substantial reduction in time step size \cite{Class02}. This feature is not desirable for simulations over a long period of time, usually encountered in groundwater remediation or transport of nuclides in a nuclear waste repository. In addition to flash calculations and PVS, there are other formulations to handle phase transitions such as negative saturation \cite{Abadpour08}, and introduction of persistent primary variables \cite{Bourgeat08,Marchand12,Neumann13}. 

Recently, a new approach has been developed for handling the phase transitions by formulating the system of equations as a nonlinear complementarity problem (NCP) \cite{BenGharbia14,Lauser11,Marchand14}. In contrast to PVS, NCP has the advantage that the set of primary variables is consistent throughout the simulation, and no primary variable switching is needed. NCP requires a complementary function, referred to as a C-function, employed to rewrite the inequality constraints for the thermodynamic state as a non-smooth equation, which requires a \textit{semi-smooth} Newton method \cite{Aganagi84,Qi93,QiSun93} to solve. Most of the previous work in multiphase flow using the NCP approach employs the minimum function as the C-function due to its simplicity for implementation and the fact that it is piecewise linear with respect to the arguments. Even though the semi-smooth Newton method applied to the NCP using the minimum function as a C-function is observed to have quadratic convergence for simple problems in porous media (see \cite{BenGharbia14}), we find that it exhibits poor convergence and even diverges for standard benchmark problems, as well as examples considered in this work. An alternative to the minimum function is the Fischer-Burmeister function, which has recently been employed as the C-function for NCP formulation of incompressible two-phase flow in \cite{Yang17}. As we will show, this choice of C-function can help mitigate the lack of robustness observed in using the Newton-min algorithm for NCP formulation of compositional two-phase flow with phase transitions. We then draw on this experience and develop a new method for the nonlinear solve based on a smooth version of the Fischer-Burmeister function. Our method can be considered a variant of the Jacobian smoothing method summarized in \cite{Facchinei03}. Compared to the non-smooth approaches that use the minimum and the Fischer-Burmeister functions, our new method is more robust and efficient for problems with highly heterogeneous media, and it also scales optimally with problem size.

We consider a two-phase, two-component system with phase transitions as our model problem. We describe this model in detail\ in \cref{sec:statement}, and in \cref{sec:ncp}, we describe the NCP formulation for it. We briefly review the semi-smooth Newton framework and introduce our new algorithm in \cref{sec:algorithm}. In \cref{sec:results}, several numerical tests are presented that demonstrate the robustness and scalability of the new algorithm. Some concluding remarks as well as future work are presented in \cref{sec:conclusions}.

\section{Problem Statement}
\label{sec:statement}
\subsection{Governing Equations}
%\subsection{Simplified 2-phase, 2-component model}
In this work, we consider a simplified two-phase two-component model with phase transitions. The phases consist of liquid and gas, and the components are water and hydrogen. We also make the following assumptions: (1) water does not vaporize so the gas phase contains only hydrogen, and (2) the amount of hydrogen dissolved in the liquid phase is small. For the two components, the mass conservation equations read
\begin{align}
&\Porosity \dfrac{\partial (\DensityComp{l}{w} \Sat{l})}{\partial t} + \nabla \cdot (\DensityComp{l}{w} \DarcyV{l} - \DiffusionFlux{l}{h}) = 0 \label{eq:water_mass_conservation},\\
&\Porosity \dfrac{\partial (\DensityComp{l}{h} \Sat{l} + \DensityComp{g}{h} \Sat{g})}{\partial t} + \nabla \cdot (\DensityComp{l}{h} \DarcyV{l} + \DensityComp{g}{h} \DarcyV{g} + \DiffusionFlux{l}{h}) = 0, \label{eq:hydrogen_mass_conservation}
\end{align}
where the subscripts $l,g$ denote the liquid and gas phases, and the superscripts $w,h$ denote the water and hydrogen components, respectively. The porosity of the medium is denoted $\Porosity$, $S_\alpha, \bm{q}_\alpha$ are the saturation and velocity of phase $\alpha$, respectively; $\DensityComp{l}{h}$ is the dissolved hydrogen mass concentration in the liquid phase; and $\DiffusionFlux{l}{h}$ is the diffusion flux of hydrogen in the liquid phase. The Darcy velocity $\bm{q}_\alpha$ follows the Darcy-Muskat law
\begin{align}
\DarcyV{\alpha} = -\AbsPerm{} \Mobility{\alpha} \nabla (\Pres{\alpha} - \Density{\alpha} \Gravity{}), \hspace{5mm} \alpha = l,g ,
\end{align}
where $\AbsPerm{}$ is the absolute permeability, $\Mobility{\alpha}$, $\Pres{\alpha}$, and $\Density{\alpha}$ are the mobility, pressure, and density of phase $\alpha$, and $\Gravity{}$ is the gravitational acceleration. The mobility $\Mobility{\alpha}$ of phase $\alpha$ is defined as the ratio between the phase relative permeability $\RelPerm{\alpha}$ and the phase viscosity $\Viscosity{\alpha}$: $\Mobility{\alpha} = \RelPerm{\alpha}/\Viscosity{\alpha}$. Using Fick's law, the diffusion flux of hydrogen in liquid $\DiffusionFlux{l}{h}$ in \cref{eq:water_mass_conservation,eq:hydrogen_mass_conservation} can be expressed as
\begin{align}
\DiffusionFlux{l}{h} = -\Porosity \Sat{l} \MolecularDiff{l}{h} \nabla \DensityComp{l}{h},
\end{align}
where $\MolecularDiff{l}{h}$ is the hydrogen molecular diffusion coefficient in liquid. Since we assume incompressibility of the liquid phase, the mass density of the water component in the liquid phase is constant, i.e. $\DensityComp{l}{w} = \DensityComp{w}{std}$. To capture capillarity effects, the jump in the pressure at the interface of the two phases is modeled by the relation
\begin{align}
&\Pres{g} = \Pres{l} + \Pres{c}(\Sat{l})
\end{align} 
where $\Pres{c}$ is the capillary pressure. Additionally, we have the constraints 
\begin{align}
\Sat{l} + \Sat{g} = 1. \label{eq:sat_constraint}
\end{align}
To close the model, we also need a set of equations for the thermodynamic equilibrium when the gas phase is present, i.e. how much hydrogen can dissolve into the liquid phase at a certain pressure. Assuming low solubility of hydrogen in the liquid phase, Henry's law can be used to connect the gas pressure $\Pres{g}$ and the dissolved hydrogen mass concentration in liquid $\DensityComp{l}{h}$:
\begin{align}
\DensityComp{l}{h} = C_h \Pres{g}, \label{eq:henrys_law}
\end{align} 
where $C_h = \HenryConstant{} \MolarMass{h} = \DensityComp{w}{std} \MolarMass{h}/(\MolarMass{w}K^h)$, $\HenryConstant{}$ is the Henry's law constant, $K^h$ is a constant specific to the mixture, and $\MolarMass{i}, \; i \in \{w,h\}$, is the molar mass of the $i$-th component. Since we neglect water vapor, we can apply the ideal gas law for the gas phase. This leads to the relation
\begin{align}
\DensityComp{g}{h} = \Density{g} = C_v \Pres{g},
\end{align}
where $C_v$ is a constant and $C_v = \MolarMass{h}/(\IdealGasConstant{}\Temp{})$; $\Temp{}$ is the temperature and $\IdealGasConstant{}$ the ideal gas constant.

\subsection{Relative Permeabilities and Capillary Pressure}
We employ the nonlinear Van Genuchten \cite{VanGenuchten80} model for relative permeabilities and capillary pressure:
\begin{align}
&\RelPerm{l} = \sqrt{\Sat{le}}\Big(1 - \big(1 - \Sat{le}^{1/m}\big)^m\Big)^2, \hspace{5mm} \RelPerm{g} = \sqrt{1 - \Sat{le}}\Big( 1- \Sat{le}^{1/m} \Big)^{2m}, \\
&\Pres{c} = \Pres{r}\Big(\Sat{le}^{-1/m} - 1\Big)^{1/n}, \\
&\Sat{le} = \dfrac{1 - \Sat{l}}{1 - \Sat{lr} - \Sat{gr}}, \hspace{5mm} m = 1 - \dfrac{1}{n},
\end{align}
where $\Pres{r}$ is the entry pressure. Notice that the function $\Pres{c}(\Sat{l})$ in the Van Genuchten model is only defined for $\Sat{l} \in [\Sat{lr},1-\Sat{gr}]$ and $\Pres{c}^\prime$ is unbounded near $\Sat{lr}$ and $1-\Sat{gr}$. Thus, it is necessary to modify the model to limit the growth of $\Pres{c}^\prime$ and extend it for $\Sat{l} \in \mathbb{R}$, since the value of $\Sat{l}$ can become larger than $1-\Sat{gr}$ or less than $\Sat{lr}$ during the nonlinear iteration. We use the following regularization as presented in \cite{Marchand14} with parameter $\epsilon = 10^{-5}$:
%\begin{itemize}
%\item For $S_g \in [S_{gr}, 1 - S_{lr}]$
%\begin{align}
%& \tilde{S} := S_{gr} + (1 - \epsilon)(S_g - S_{gr}) + \dfrac{\epsilon}{2}(1 - S_{gr} - S_{lr}) \\
%& \tilde{P}_c(S_g) = P_c(\tilde{S}) - P_c\Big( S_{gr} + \dfrac{\epsilon}{2}(1 - S_{gr} - S_{lr}) \Big)
%\end{align}
%\item For $S_g < S_{gr}$
%\begin{align}
%\tilde{P}_c(S_g) = \tilde{P}_c(S_{gr}) + \tilde{P}_c^\prime (S_{gr})(S_g - S_{gr})
%\end{align}
%\item For $S_g > 1 - S_{lr}$
%\begin{align}
%\tilde{P}_c(S_g) = \tilde{P}_c(1 - S_{lr}) + \tilde{P}_c^\prime (1 - S_{lr})(S_g - 1 + S_{lr})
%\end{align}
%\end{itemize}
In this regularization, for the saturation that is outside of the domain, capillary pressure is computed by a linear extrapolation from the regularization points $\Sat{lr} + \mathcal{O}(\epsilon)$ and $1 - \Sat{gr} - \mathcal{O}(\epsilon)$ with the slopes $\Pres{c}^\prime (\Sat{lr} + \mathcal{O}(\epsilon))$ and $\Pres{c}^\prime (1 - \Sat{gr} - \mathcal{O}(\epsilon))$, respectively.

\subsection{Primary Variables}\label{subsec:primary_variables}
There are many ways to choose a set of primary variables, depending on the problem formulation and applications. In our model example, a convenient choice is the liquid pressure, liquid saturation, and the concentration of hydrogen in the liquid phase. We then have our solution vector $u = \{\Pres{l}, \Sat{l}, \DensityComp{l}{h}\}$. Unlike in other methods such as primary variable switching, for NCP, the choice of primary variables is fixed throughout the simulation.

\section{Nonlinear Complementarity Problem}
\label{sec:ncp}
In its simplest form, a nonlinear complementarity problem with respect to a smooth function $f: \mathbb{R}^N \mapsto \mathbb{R}^N$ is to find a vector $u \in \mathbb{R}^N$ such that
\begin{align}
u \ge 0, \hspace{5mm}
f(u) \ge 0, \hspace{5mm}
u^T f(u) = 0, \label{eq:complementarity_equation}
\end{align}
A slightly more general form of the last equation in \cref{eq:complementarity_equation} reads
\begin{align}
g(u)^T f(u) = 0, \label{eq:general_ncp}
\end{align}
where $g: \mathbb{R}^N \rightarrow \mathbb{R}^N$ is another smooth function. As we have mentioned in \cref{sec:statement}, for the solution of \cref{eq:water_mass_conservation,eq:hydrogen_mass_conservation} to be valid, the pressure, saturation, and hydrogen concentration in the liquid phase must satisfy the constraints in  \cref{eq:sat_constraint,eq:henrys_law}. These conditions can be reformulated as an NCP as follows:
\begin{align}
1 - \Sat{l} \geq 0, \hspace{3mm} C_h \Pres{g} - \DensityComp{l}{h} \geq 0, \hspace{3mm} (1 - \Sat{l})(C_h \Pres{g} - \DensityComp{l}{h}) = 0, \label{eq:ncp_problem}
\end{align}
where $u = \begin{pmatrix}
\Pres{l} & \Sat{l} & \DensityComp{l}{h}
\end{pmatrix}^T$, and the functions in \cref{eq:general_ncp} are $g(u) = 1-\Sat{l}$, $f(u) = C_h \Pres{g} - \DensityComp{l}{h}$. A very popular approach to solve \cref{eq:ncp_problem} is to transform it into a semi-smooth nonlinear equation via a complementarity function (also called C-function) $\Phi (a,b): \mathbb{R}^2 \rightarrow \mathbb{R}$, which satisfies
\begin{align}
\Phi (a,b) = 0 \iff a \geq 0, \hspace{3mm} b \geq 0, \hspace{3mm} ab = 0. \label{eq:cfun}
\end{align}
We can extend the definition of $\Phi$ \cref{eq:cfun} from $\mathbb{R}^2$ to $\mathbb{R}^N$, where $N$ is also the number of elements in the mesh, by applying it \cref{eq:cfun} componentwise to $1 - \Sat{l}$ and $C_h \Pres{g} - \DensityComp{l}{h}$ and obtain the nonlinear system
\begin{align}
\Theta (u) = \begin{pmatrix}
\Phi (1-(\Sat{l})_1, (C_hP_g - \DensityComp{l}{h})_1) \\
\Phi (1-(\Sat{l})_2, (C_hP_g - \DensityComp{l}{h})_2) \\
\cdots \\
\Phi (1-(\Sat{l})_N, (C_hP_g - \DensityComp{l}{h})_N)
\end{pmatrix}. \label{eq:discrete_cfun}
\end{align}
Then, solving the NCP problem in \cref{eq:ncp_problem} is equivalent to solving $\Theta (u) = 0$. As discussed in the next section, combining this complementarity condition with the discrete PDEs enables use of a nonlinear solution algorithm that automatically enforces the constraints \eqref{eq:sat_constraint}-\eqref{eq:henrys_law}.
%with functions $a(u) = 1 - \Sat{l}$ and $b(u) = C_h \Pres{g} - \DensityComp{l}{h}$. 
There are many examples of C-functions \cite{Facchinei03}. In this work, we focus on two popular choices
\begin{align}
&\Phi _{\min} (a,b) = \min(a,b) \label{eq:min_ncp}\\
&\Phi _{FB} (a,b) = \sqrt{a^2 + b^2} - (a+b) \hspace{3mm}\text{(Fischer-Burmeister)}\label{eq:fb_ncp}
\end{align}
%%If we further introduce a merit function \footnote{A merit function for $H(x) = 0$ is any function $M:\mathbb{R}^n \mapsto \mathbb{R}$ such that $M(x) \geq 0$ for any $x \in \mathbb{R}^n$ and $M(x) = 0$ iff $H(x) = 0$.}
%%\begin{align}
%%\psi (a,b) = \dfrac{1}{2}\norm{\Porosity (a,b)}^2,
%\end{align}
%then finding the solution to \cref{eq:ncp_problem} is also equivalent to minimizing $\psi$. 
The minimum function is convenient because it is piecewise linear with respect to the variables $a$ and $b$, which simplifies the computation of the Jacobian in each nonlinear iteration \cite{Lauser11}. When the gas phase is not present, \cref{eq:min_ncp} reduces to $1 - \Sat{l} = 0$. When the gas phase appears, $1 - \Sat{l} > 0$ and the constraint equation is governed by Henry's law \cref{eq:henrys_law}. Recently, the Fischer-Burmeister function has been shown to have good performance for the case of incompressible two-phase flow \cite{Yang17}.

\par Although we do not employ any line search strategy in this work, we note that compared to the Fischer-Burmeister function, the minimum function is less useful with respect to globalization with line search strategies \cite{Facchinei03}. As shown and discussed in \cite{BenGharbia12,BenGharbia13,BenGharbia18}, global semi-smooth Newton methods may diverge even for linear C-functions if the starting point is not close enough to a solution. In global semi-smooth Newton approaches, this can be handled using line search, in which a \textit{merit function} is used to enhance robustness by evaluating the quality of new iterates. In this scenario, it is desirable for the natural merit function $\Psi = \norm{\Phi}^2$ to be smooth because the search direction at each nonlinear iteration is usually chosen based on the derivative of $\Psi$ \cite{Kelley18}. The merit function associated with the minimum function $\Psi _{\min} = \norm{\Phi _{\min} (a,b)}^2$, however, does not satisfy this condition. In contrast, the Fischer-Burmeister merit function $\Psi _{FB} = \norm{\Phi _{FB} (a,b)}^2$ is continuously differentiable as observed in \cite{DeLuca96}.
%\footnote{A merit function for $H(x) = 0$ is any function $M:\mathbb{R}^n \mapsto \mathbb{R}$ such that $M(x) \geq 0$ for any $x \in \mathbb{R}^n$ and $M(x) = 0$ iff $H(x) = 0$.}. 
%Since the step size $\alpha_k$ is usually chosen based on the derivative of the merit function, it is desirable for the merit function to be continuously differentiable. 

%\footnote{A merit function for $H(x) = 0$ is any function $M:\mathbb{R}^n \mapsto \mathbb{R}$ such that $M(x) \geq 0$ for any $x \in \mathbb{R}^n$ and $M(x) = 0$ iff $H(x) = 0$.}. In constrained optimization, it is possible that the value of the cost function $f(x)$ is decreased at the iterate $x_{k+1}$ but at the cost of an increase in violation of the constraints. Thus, a merit function is employed as a measure of progress to evaluate the choice of the step size $\alpha_k$ in the update $x_{k+1} = x_k + \alpha_k\bm{p}_k$ where $\bm{p}_k$ is a search direction. 

\section{Solution Algorithm}
\label{sec:algorithm}
We consider solving the coupled system consisting of \cref{eq:water_mass_conservation,eq:hydrogen_mass_conservation,eq:ncp_problem} fully implicitly. We use a cell-centered finite volume method for spatial discretization, as it is a natural way to preserve the mass conservation property of the balance \cref{eq:water_mass_conservation,eq:hydrogen_mass_conservation}. In addition, it can deal with the case of discontinuous permeability coefficients, and it is relatively straightforward to implement. 
%Under appropriate assumptions, this method also falls into the mixed finite element framework \cite{Peaceman77,Russell83}. 
For the time domain, we employ the backward Euler method to avoid a CFL stability restriction on the time step. Because this method is unconditionally stable, it also allows us to experiment with variable time stepping, which can significantly reduce execution time.
\subsection{Semi-smooth Newton Method}
We want to solve the system $R(u) = 0$ where $R(u)$ is the residual function given by
\begin{align}
R(u) = \begin{dcases}
H(u) & \text{(from the PDEs)} \\
\Theta(u) & \text{(from the constraints)}  \\
\end{dcases} \label{eq:ncp_system}
\end{align}
This is a system of $3N$ equations and $3N$ unknowns, where $N$ is the number of elements in the mesh. A standard approach for solving nonlinear systems of equations is Newton's method, which requires solution of a linear system at each iteration $k$:
\begin{align}
\dfrac{\partial R}{\partial u}\Big|_{u=u_k} \delta u = - R(u_k) .
\end{align}
This method requires that the Jacobian $\partial R/\partial u$ be defined everywhere. In the NCP formulation, the constraints $\Theta$ are not differentiable when there is phase transition as the solution changes from satisfying $1-\Sat{l} = 0$ to $C_hP_g - \DensityComp{l}{h} = 0$. To address this, we will consider a semi-smooth Newton method, which is similar to Newton's method, except the derivative $\Theta^\prime$ is replaced by a member of the \textit{subdifferential} $\partial \Theta$ when $\Theta$ is not differentiable. Let $F: \mathbb{R}^n \mapsto \mathbb{R}^n$ be a locally Lipschitz-continuous function and $D_F$ be the set where $F$ is differentiable; the B-subdifferential of $F$ at $x$ is defined as the set
\begin{align*}
\partial _B F(x) := \{G \in \mathbb{R}^{n\times n}: \exists \; x_k \in D_F \text{ with } x_k \rightarrow x, \nabla F(x_k) \rightarrow G \}\;.
\end{align*}
Below is the algorithm for the general semi-smooth Newton method (see \cite{Facchinei03}).
\begin{algorithm}[H]
\While{k $<$ max\_iter and res $>$ tol}{
  (1) Given $u^0$, $k = 0$ \\
  (2) Select an element $J_k \in \partial _B \Theta(u^k)$\\
  (3) Solve the system
     \begin{flalign*}
       \begin{pmatrix}
       H^\prime(u^k) \\
       J_k (u^k)
       \end{pmatrix} \bigtriangleup u^k &= 
       \begin{pmatrix}
       - H(u^k)\\
       - \Theta (u^k)
       \end{pmatrix}&
     \end{flalign*}
  (4) Update $u^{k+1}$: $u^{k+1} = u^k + \bigtriangleup u^k$
}
\caption{General semi-Smooth Newton method.}\label{algo:sn_fb}
\end{algorithm}
To compute $J_k$ in the algorithm above, one can use an active set strategy \cite{Hintermuller02}. For multiphase flow with phase appearance and disappearance, the idea is to define the set of indices for the cells in which the gas phase is present (see \cite{BenGharbia14,Lauser11}). Let $A^k := \lbrace j: 1-(\Sat{l})_j \geq (C_h \Pres{g} - \DensityComp{l}{h})_j \rbrace$,  $I^k := \lbrace j: 1-(\Sat{l})_j < (C_h \Pres{g} - \DensityComp{l}{h})_j \rbrace$. Then for the minimum function, the $j$th row of $J_k$ is equal to 
\begin{align}
\begin{dcases}
\dfrac{\partial}{\partial u} a(u)_j \hspace{3mm} \text{if} \hspace{3mm} j \in I^k \\
\dfrac{\partial}{\partial u} b(u)_j \hspace{3mm} \text{if} \hspace{3mm} j \in A^k \\
\end{dcases}
\end{align}
where $a(u)_j = 1 - (\Sat{l})_j$ and $b(u)_j = (C_h \Pres{g} - \DensityComp{l}{h})_j$. For the Fischer-Burmeister function, we can compute the $j$th row of $J_k$  as follows:
\begin{align}
\begin{dcases}
\dfrac{1}{\sqrt{a(u)_j^2 + b(u)_j^2}}\Big(a(u)_j \dfrac{\partial}{\partial u} a(u)_j + b(u)_j \dfrac{\partial}{\partial u} b(u)_j  \Big) - 
\Big( \dfrac{\partial}{\partial u} a(u)_j + \dfrac{\partial}{\partial u} b(u)_j \Big) \hspace{3mm} &\text{if} \hspace{1mm} a(u)_j^2 + b(u)_j^2 \neq 0 \\
(\alpha_i - 1)\dfrac{\partial}{\partial u} a(u)_j + (\beta_i - 1)\dfrac{\partial}{\partial u} b(u)_j \hspace{3mm} &\text{otherwise}
\end{dcases} \label{eq:subdiff_fb}
\end{align}
where for all $i$ such that $a(u)_j^2 + b(u)_j^2 = 0$, $\alpha_i$ and $\beta_i$ are arbitrary nonnegative constants satisfying $\alpha_i^2 + \beta_i^2 = 1$. For a more complete treatment of semi-smooth Newton methods, we refer to \cite{Facchinei03}.

\subsection{Jacobian Smoothing Method}
An alternative to the semi-smooth approach is to employ a smooth approximation to the non-smooth function $\Theta$. The idea was originally developed in the context of variational inequalities \cite{Chen98, Chen99} and generalized to more general nonlinearities and infinite dimensions in \cite{Chen01}. Let $G: \mathbb{R}^n \times \mathbb{R}_{+} \mapsto \mathbb{R}^n$ such that for any $\tau > 0$, $G(\cdot,\tau)$ is continuously differentiable on $\mathbb{R}^n$ and 
\begin{align}
\norm{\Theta(u) - G(u,\tau)} \rightarrow 0, \hspace{3mm} \text{as} \hspace{3mm} \tau \rightarrow 0.
\end{align} 
Then, given a sequence $\tau ^k$, $k = 0,1,2,...$, we can solve the system in \cref{eq:ncp_system} inexactly using $G^\prime(u^k,\tau ^k)$ as an approximation to the generalized Jacobian $J_k = \partial _B \Theta(u^k)$. In this work, we explore a smooth approximation to the Fischer-Burmeister functions given by
\begin{align}
%&G_{\min} (u,\tau) = \dfrac{1}{2}(-\sqrt{(a-b)^2 + 2\tau} + a + b)\\
G_{FB} (u,\tau) = \sqrt{a^2 + b^2 + 2\tau} - (a+b)
\end{align}
The complete algorithm is as follows:
\begin{algorithm}[H]
\While{k $<$ max\_iter and res $>$ tol}{
  (1) Given $u^0$, $k = 0$, and $\tau ^0$\\
  (2) Solve the system
     \begin{flalign*}
       \begin{pmatrix}
       H^\prime(u^k) \\
       G^\prime(u^k,\tau^k)
       \end{pmatrix} \bigtriangleup u^k &= 
       \begin{pmatrix}
       - H(u^k)\\
       - \Theta (u^k)
       \end{pmatrix}&
     \end{flalign*}
  \hspace{-2mm}(3) Update the smoothing parameter $\tau$\\
     \hspace{6mm} $\tau^{k+1} = \beta \tau^k$ for $\beta \in (0,1)$\\
  (4) Update $u^{k+1}$ \\
     \hspace{6mm}$u^{k+1} = u^k + \bigtriangleup u^k$
}
\captionsetup[algorithm]{labelsep=colon}
\caption{Jacobian Smoothing Method.}\label{algo:jsm}
\end{algorithm}
There also exist smooth approximations to the minimum function. In particular, we experimented with the Chen-Harker-Kanzow-Smale smoothing \cite{Chen93,Kanzow96},
\begin{align}
G_{min} (u,\tau) = (a+b) - \sqrt{(a-b)^2 + 4\tau}.
\end{align}
However, in our experience, this smooth version of the minimum function does not improve the convergence of the semi-smooth Newton's method significantly, and we do not include the results here.

\subsection{Linear System}
Assuming that each physical variable is ordered lexicographically, each step of the nonlinear iteration (step 3 in \cref{algo:sn_fb} and step 2 in \cref{algo:jsm}) requires solution of a large sparse, non-symmetric, indefinite linear system of the form
\begin{align}
\begin{pmatrix}
A_{11} & A_{12} & A_{13} \\
A_{21} & A_{22} & A_{23} \\
A_{31} & A_{32} & A_{33}
\end{pmatrix} 
\begin{pmatrix}
u_1 \\
u_2 \\
u_3
\end{pmatrix} = \begin{pmatrix}
f_1 \\
f_2 \\
f_3
\end{pmatrix}. \label{eq:linear_system}
\end{align}
where $u_1, u_2, u_3$ are the corrections to $\Sat{l}, \Pres{l}, \DensityComp{l}{h}$, respectively. The matrices in the first two rows are the discretized version of the linearized operators from the PDEs, and the last row corresponds to the discrete derivative of the complementarity constraint equation introduced in \cref{eq:discrete_cfun}. Iterative methods such as GMRES \cite{Saad86} are the only viable option to solve the system above, and preconditioning is critical for fast convergence. Here, all of our experiments use GMRES preconditioned with hypreMGR (see \cite{Bui18}), an AMG solver and preconditioner based on multigrid reduction and designed for systems of PDEs. Unlike ILU preconditioners in which one only needs to specify the level of fill (in ILU($k$)) or the threshold tolerance (in ILU($t$)), hypreMGR requires extra information regarding the block structure of the system and the order of reduction.
\par There exists a small but important difference in the structure of the matrices $A$ created using the Jacobian smoothing method and the semi-smooth Newton methods. For the semi-smooth Newton methods with an active set strategy, the diagonal of the block $A_{33}$ contains zeros for the cells that are devoid of the gas phase. In contrast, for $\tau > 0$, the diagonal of the block $A_{33}$ is guaranteed to be nonzero for the Jacobian smoothing method regardless of the existence of phase transitions. Thus, the Jacobian smoothing method requires one fewer reduction step for hypreMGR, compared to semi-smooth Newton approaches, which leads to a decrease in both the number of GMRES iterations and execution time, as will become evident from the results presented in \cref{subsec:scaling}.

\section{Numerical Results}
\label{sec:results}
In this section, we describe the results of numerical experiments for solving the NCP systems using both the semi-smooth Newton's approach using the minimum and the Fischer-Burmeister functions, and the Jacobian smoothing method with the smooth Fischer-Burmeister function for the NCP formulation. All of these methods are implemented in Amanzi, a parallel open-source multi-physics C++ code developed as a part of the ASCEM project \cite{Amanzi}. Although Amanzi was first designed for simulation of subsurface flow and reactive transport, its modular framework and concept of process kernels \cite{Coon14} allow new physics to be added relatively easily for other applications. The compositional two-phase flow simulator employed in this work is one such example. Amanzi works on a variety of platforms, from laptops to supercomputers. It also leverages several popular packages for mesh infrastructure and solvers through a unified input file. GMRES is also provided within Amanzi while hypreMGR is employed through HYPRE.

\par This section has three parts. In the first part, we show results using two benchmark problems that demonstrate the effectiveness of the NCP approach in handling phase appearance and disappearance. In the second part, we report the results for both two- and three-dimensional cases with highly heterogeneous media. In the last part, we perform a scalability study of the algorithms. Parallel test cases are run on Syrah, a Cray system with 5,184 Intel Xeon E5-2670 cores at the Lawrence Livermore National Laboratory Computing Center. Amanzi and other libraries are compiled with OpenMPI 1.6.5 and gcc-4.9.2.

\par For all the simulations presented here, the convergence tolerance for the nonlinear solve is $||F(x)|| \le 10^{-6}$, and the linear tolerance for GMRES is $||J\delta u_k - F(u_k)|| \le 10^{-12}||F(u_k)||$, which is the default in Amanzi. Depending on the performance of the nonlinear solver, a heuristic for choosing the time step is used: if the number of nonlinear steps (NS) required at a given time are less than 10, then the next time step $dt_{next}$ is doubled, $dt_{next} = 2\cdot dt$; if NS $\in [11,15]$, then the time step is kept fixed, $dt_{next} = dt$; and if NS is greater than 15, then the time step is halved $dt_{next} = dt / 2$. The maximum number of nonlinear iterations is $max\_iter = 20$.

\subsection{Benchmark problems}\label{subsec:benchmark_problems}
\begin{figure}
\centering
\includegraphics[width=0.7\textwidth]{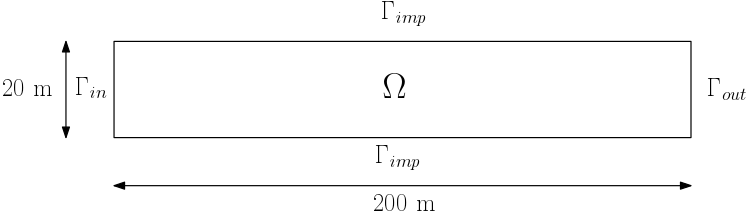}
\captionsetup{justification=centering}
\caption{Core domain for the gas infiltration example.} \label{fig:radon1}
\end{figure}

These tests are derived from the MoMaS benchmark project \cite{BourgeatBenchmark13}, which is designed to evaluate the effectiveness of different approaches for handling gas phase appearance and disappearance. Pure hydrogen is injected into a two-dimensional homogeneous porous domain $\Omega$, which was initially 100\% saturated with pure water. The domain is a rectangle of size $200\text{m} \times 20\text{m}$, and it is discretized only in the horizontal direction, leading to a quasi one-dimensional problem. There are three types of boundaries : $\Gamma _{\text{in}}$ on the left side is the inflow boundary; $\Gamma _{\text{out}}$ on the right side is the outflow boundary; and $\Gamma _{\text{imp}}$ at the top and bottom is the impervious boundary (see \cref{fig:radon1}).
\noindent There are no source terms inside the domain, and denoting $\Flux{w} = \DensityComp{l}{w} \AbsPerm{} \Mobility{l} \nabla \Pres{l} - \DiffusionFlux{l}{h}$ and $\Flux{h} = \DensityComp{l}{h} \AbsPerm{} \Mobility{l} \nabla \Pres{l} + \DensityComp{g}{h} \AbsPerm{} \Mobility{g} \nabla \Pres{g} + \DiffusionFlux{l}{h}$, the boundary conditions are as follows
\begin{itemize}
\item No flux on $\Gamma _\text{imp}$
\begin{align}
\Flux{w} \cdot \nu = 0 \text{   and   } \Flux{h} \cdot \nu = 0
\end{align}
\item Injection of hydrogen on the inlet $\Gamma _{\text{in}}$
\begin{align}
\Flux{w} \cdot \nu = 0 \text{   and   } \Flux{h} \cdot \nu = 5.57 \times 10^{-6} \hspace{2mm} kg/m^2/year
\end{align}
\item Fixed liquid saturation and pressure on the outlet
\begin{align}
\Pres{l} = 10^6 \hspace{2mm} Pa, \hspace{3mm} \Sat{l} = 1, \hspace{3mm} \DensityComp{l}{h} = 0
\end{align}
\end{itemize}
Initial conditions are uniform throughout the domain, corresponding to a stationary state of saturated liquid and no hydrogen injection,
\begin{align}
\Pres{l} = 10^6 \hspace{2mm} Pa, \hspace{3mm} \Sat{l} = 1, \hspace{3mm} \DensityComp{l}{h} = 0.
\end{align}
The values of the physical parameters are given in \cref{tab:param_list_radon1}.
\begin{table}[ht]
\begin{minipage}[b]{0.35\linewidth}
\centering
{\renewcommand{\arraystretch}{1.2}
\begin{tabularx}{\textwidth}{ l X }
  \hline
  $\AbsPerm{}$ & $5 \times 10^{-20}$ $m^2$ \\
  $\Porosity$ & 0.15 \\
  $\MolecularDiff{l}{h}$ & $3 \times 10^{-9}$ $m^2/s$ \\
  $\Viscosity{l}$ & $1 \times 10^{-9}$ Pa s\\
  $\Viscosity{g}$ & $9 \times 10^{-6}$ Pa s\\
  $\HenryConstant{}$ & $7.65 \times 10^{-6}$ mol/Pa/$m^3$ \\
  $\MolarMass{h}$ & $2 \times 10^{-3}$ kg/mol\\
  $\MolarMass{w}$ & $1 \times 10^{-2}$ kg/mol\\
  $\DensityComp{l}{w}$ & $10^3$ kg/$m^3$ \\
  \hline
\end{tabularx}}
\caption{Parameter Values} \label{tab:param_list_radon1}
\end{minipage}\hfill
\begin{minipage}[b]{0.6\linewidth}
\centering
  \includegraphics[width=0.9\textwidth]{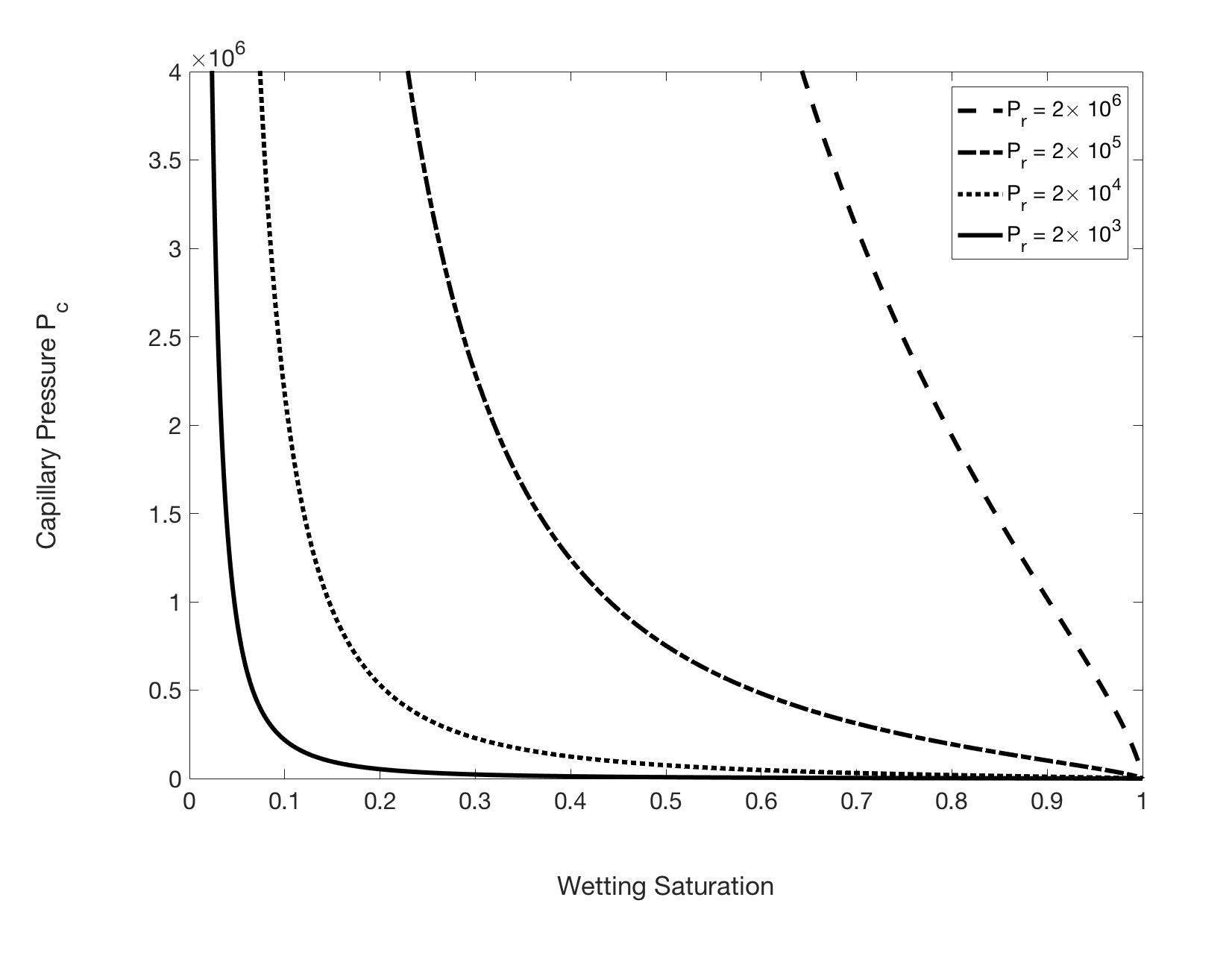}
  \captionsetup{justification=centering}
  \captionof{figure}{Capillary pressure curves for different entry pressure $\Pres{r}$.} \label{fig:capillary_curves}
\end{minipage}
\end{table}
\par
\begin{figure}
\centering
\includegraphics[width=0.5\textwidth]{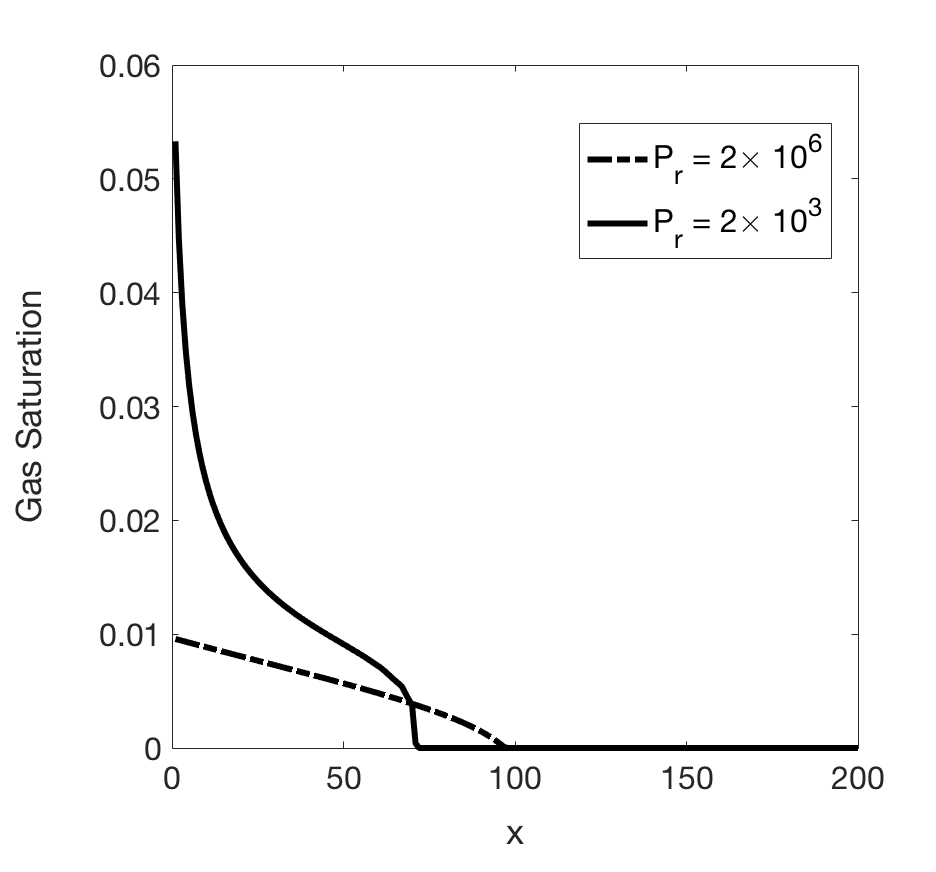}
\captionsetup{justification=centering}
\caption{Gas infiltration into the domain for two different capillary pressure curves after 100,000 years.} \label{fig:gas_profiles}
\end{figure}

 \Cref{fig:capillary_curves} shows the Van Genuchten capillary pressure curve for different values of the entry pressure $\Pres{r}$. These parameter values, along with others in the Van Genuchten model, depend on the porous material. For example, the MoMaS benchmark problem (test case 1 in \cite{BourgeatBenchmark13}) uses $\Pres{r} = 2\cdot 10^6$ for a very dense rock with extremely low permeability of $\AbsPerm{} = 5\cdot 10^{-20}$. In other applications including CO$_2$ sequestration and reservoir simulation, the material is much more permeable and $\Pres{r} = 2\cdot 10^3$ would produce the capillary pressure curves typically used (see \cite{Class09,Nordbotten12}). Other parameters for the Van Genuchten model are $\Sat{lr} = 0.4$, $\Sat{gr} = 0$, and $n = 1.49$. The effect of capillary pressure on the solution is shown in \cref{fig:gas_profiles}, in which the gas saturation throughout the domain is plotted at 100,000 years for $\Pres{r} = 2\cdot 10^6$ and $\Pres{r} = 2\cdot 10^3$. The smaller $\Pres{r}$ is, the steeper the curve becomes near $\Sat{l} = 0$, and that also makes the problem more difficult to solve. For the MoMaS benchmark case with $\Pres{r} = 2\cdot 10^6$, the gas saturation curve exhibits a more gradual transition from the unsaturated to the saturated region. In contrast, for the difficult case of $\Pres{r} = 2\cdot 10^3$ Pa, the gas saturation changes very quickly both at the injection point and at the interface with the saturated region. We note that the simulation results in \cref{fig:gas_profiles} match well with those in \cite{Cuveland15,Marchand14}. A comparison of the performance of the three solution methods is shown in \cref{tab:result_momas,tab:result_neumann}. TS, NS denote the total number of successful time steps and nonlinear iterations, respectively, and the numbers in parentheses are the number of failed time steps and nonlinear iterations. Failed time steps are those in which the method diverges or does not converge within the allowed maximum number of iterations; failed nonlinear iterations correspond to the iterations spent during the failed time steps. For both of these benchmark problems, an initial smoothing parameter $\tau = 10^{-6}$ and a reduction ratio $\beta = 0.1$ are used for the smooth Fischer-Burmeister approach.
%\par For the MoMaS gas injection benchmark problem with $\Pres{r} = 2\cdot 10^6$ Pa, the results in \cref{tab:result_momas} show that for the nonlinear solve, the Fischer-Burmeister function (without smoothing) is a better choice for the C-function than the minimum function. It allows for larger time steps, indicated by the smaller numbers of time steps needed to run the simulation to $10^5$ and $5\cdot 10^5$ years, respectively. Compared to the minimum function, the Fischer-Burmeister function also reduces the number of nonlinear iterations by about 10\% and 13\%, respectively. The smooth Fischer-Burmeister function achieves about the same performance as the Fischer-Burmeister for this problem, but we also found that the former is slightly better at the beginning of the simulation (time period of $0-10^5$ years) when the gas saturation rises very quickly and the problem is highly nonlinear. From $10^5$ to $5\cdot 10^5$ years, the gas saturation still increases, but at a much slower rate, and it takes both methods about the same number of time steps to simulate both periods.
\par For the MoMaS gas injection benchmark problem with $\Pres{r} = 2\cdot 10^6$ Pa, the results in \cref{tab:result_momas} show that for the nonlinear solve, the Fischer-Burmeister function (without smoothing) does not show any improvement over the minimum function. It registers the same numbers of time steps needed to run the simulation both to $10^5$ and to $5\cdot 10^5$ years. In contrast, the smooth Fischer-Burmeister function achieves the same performance up to $T = 10^5$ years, and it reduces both the number of time steps and nonlinear iterations by about 20\% for the full simulation. This suggests that the smooth Fischer-Burmeister function is better for simulating long time periods, when the gas phase infiltrates a larger portion of the domain.
\begin{table}
\centering
{\renewcommand{\arraystretch}{1.2}
\begin{tabular}{|c|c|c|c|c|c|c|}
\hline
  \multirow{2}{*}{End Time (years)} & \multicolumn{2}{c|}{$\min$} &   \multicolumn{2}{c|}{FB} & \multicolumn{2}{c|}{Smooth FB} \\
\cline{2-7}
  & TS & NS & TS & NS & TS & NS \\
\hline
$10^5$ & 5 (0) & 35 (0) & 5 (0) & 35 (0) & 5 (0) & 36 (0) \\
$5 \cdot 10^5$ & 10 (0) & 80 (0) & 10 (0) & 80 (0) & 8 (0) & 63 (0) \\ 
\hline
\end{tabular}}
\captionsetup{justification=centering}
\caption{Performance of the nonlinear solver for the capillary pressure model with $\Pres{r} = 2\times 10^6$ Pa with mesh size of 200.} \label{tab:result_momas}
\end{table}
\begin{table}
\centering
{\renewcommand{\arraystretch}{1.2}
\begin{tabular}{|c|c|c|c|c|c|c|}
\hline
  \multirow{2}{*}{Mesh size} & \multicolumn{2}{c|}{$\min$} &   \multicolumn{2}{c|}{FB} & \multicolumn{2}{c|}{Smooth FB} \\
\cline{2-7}
  & TS & NS & TS & NS & TS & NS \\
\hline
 200 & 37 (20) & 487 (195) & 5 (0) & 41 (0) & 5 (0) & 38 (0) \\ 
 400 & 59 (48) & 949 (440) & 6 (0) & 59 (0) & 5 (0) & 42 (0) \\
\hline
\end{tabular}}
\captionsetup{justification=centering}
\caption{Performance of the nonlinear solver for the highly nonlinear capillary pressure model with $\Pres{r} = 2\times 10^3$ Pa after 100,000 years.} \label{tab:result_neumann}
\end{table}
%\par The second example illustrates the effectiveness of the smooth Fischer-Burmeister function in handling phase transitions for highly nonlinear problems. We compare the performance of the three different strategies and show the results in \cref{tab:result_neumann}. The semi-smooth Newton method with the minimum function struggles to converge for many time steps. It requires 37 and 58 time steps in total, with 12 and 18 failed time steps for mesh sizes of $200$ and $400$, respectively. Using the Fischer-Burmeister function reduces the number of time steps by a factor of three, and it also requires only about 60\% number of nonlinear iterations. This means that on average, we can take about three times larger time step and achieve approximately 60\% faster run time with the Fischer-Burmeister function. 
The second example is exactly the same as the first example, except that we use the entry pressure $\Pres{r} = 2\cdot 10^3$ for the Van Genuchten model, which makes the problem more difficult to solve. This example the  illustrates the effectiveness of the smooth Fischer-Burmeister function in handling phase transitions for highly nonlinear problems. We compare the performance of the three different strategies and show the results in \cref{tab:result_neumann}. The semi-smooth Newton method with the minimum function struggles to converge for many time steps. It requires 37 and 59 time steps in total, with 20 and 48 failed time steps for mesh sizes of $200$ and $400$, respectively. Use of the Fischer-Burmeister function reduces the number of time steps by a factor of seven, and it also requires less than 10\% number of nonlinear iterations. This means that on average, we can take about seven times larger time step and achieve approximately 90\% decrease in run time with the Fischer-Burmeister function. 
%\par The approach using the smooth Fischer-Burmeister function registers the same number of time steps as the approach using the standard Fischer-Burmeister function and it furthers decreases the number of time steps by 10\% for the mesh size of $200$. For the larger mesh of $400$, however, the smooth Fischer-Burmeister variant shows a large improvement over the standard Fischer-Burmeister approach, requiring 36\% fewer time steps and nonlinear iterations.
\par The approach using the smooth Fischer-Burmeister function registers about the same number of time steps as the approach using the standard Fischer-Burmeister function and it furthers decreases the number of time steps by 7\% for the mesh size of $200$. For the larger mesh of $400$, however, the smooth Fischer-Burmeister variant shows a large improvement over the standard Fischer-Burmeister approach, requiring 29\% fewer nonlinear iterations.
%\begin{table}[H]
%\centering
%{\renewcommand{\arraystretch}{1.2}
%\begin{tabular}{|c|c|c|c|c|c|c|}
%\hline
%  \multirow{2}{*}{End Time (years)} & \multicolumn{2}{c|}{$\min$} &   \multicolumn{2}{c|}{FB} & \multicolumn{2}{c|}{Smooth FB} \\
%\cline{2-7}
%  & TS & NS & TS & NS & TS & NS \\
%\hline
%$10^5$ & 11 (0) & 82 (0) & 10 (0) & 74 (0) & 9 (0) & 68 (0) \\
%$5 \cdot 10^5$ & 22 (0) & 177 (0) & 19 (0) & 153 (0) & 19 (0) & 157 (0) \\ 
%\hline
%\end{tabular}}
%\captionsetup{justification=centering}
%\caption{Performance of the nonlinear solver for the capillary pressure model with $\Pres{r} = 2\times 10^6$ Pa with mesh size of 200.} \label{tab:result_momas}
%\end{table}
%\begin{table}[H]
%\centering
%{\renewcommand{\arraystretch}{1.2}
%\begin{tabular}{|c|c|c|c|c|c|c|}
%\hline
%  \multirow{2}{*}{Mesh size} & \multicolumn{2}{c|}{$\min$} &   \multicolumn{2}{c|}{FB} & \multicolumn{2}{c|}{Smooth FB} \\
%\cline{2-7}
%  & TS & NS & TS & NS & TS & NS \\
%\hline
% 200 & 37 (12) & 267 (83) & 11 (1) & 100 (8) & 11 (1) & 90 (8) \\ 
% 400 & 58 (18) & 440 (109) & 19 (3) & 182 (24) & 12 (1) & 116 (8) \\
%\hline
%\end{tabular}}
%\captionsetup{justification=centering}
%\caption{Performance of the nonlinear solver for the highly nonlinear capillary pressure model with $\Pres{r} = 2\times 10^3$ Pa after 50,000 years} \label{tab:result_neumann}
%\end{table}

\subsection{Problems with highly heterogeneous media}\label{subsec:heterogeneous}
\begin{figure}
\centering
\begin{subfigure}[t]{0.50\textwidth}
\centering
\includegraphics[width=\textwidth]{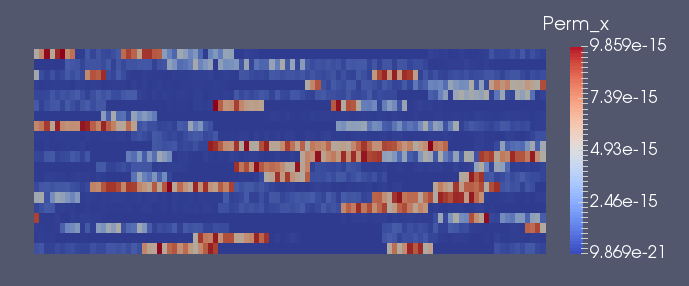}
\captionsetup{justification=centering}
\caption{Modified two-dimensional SPE10 problem.}  \label{fig:spe10_2d}
\end{subfigure}
\begin{subfigure}[t]{0.395\textwidth}
\centering
\includegraphics[width=\textwidth]{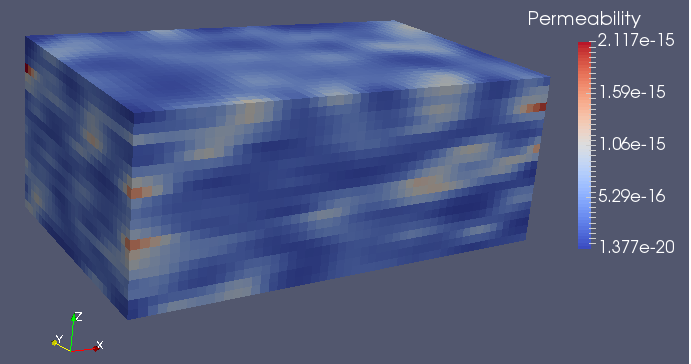}
\captionsetup{justification=centering}
\caption{Three-dimensional problem with random permeability.}  \label{fig:3d_problem}
\end{subfigure}
\caption{Heterogeneous Problems.}
\end{figure}

Here, we describe numerical experiments on two problems with highly heterogeneous permeability: (1) a modified two-dimensional SPE-10 problem, and (2) a three-dimensional problem. The permeability fields for these problems are shown in \cref{fig:spe10_2d,fig:3d_problem}. In both cases, the entry pressure for the Van Genuchten capillary pressure is chosen as $\Pres{r} = 2\times 10^3$, which corresponds to the difficult nonlinear case for the benchmark problem of \cref{subsec:benchmark_problems}.
For the first case, we modify the two-dimensional SPE10 problem \cite{Christie01} by scaling the permeability field by a constant factor of $10^{-5}$ to make the porous medium more dense. The domain is a rectangle of size $762m \times 15.24m$. Pure hydrogen is injected on the left side $\Gamma _{\text{in}} = \{0\} \times [0,15.24]$: $\Flux{w} \cdot \nu = 0$ and $\Flux{h} \cdot \nu = 5.57 \times 10^{-2} kg/m^2/year$, and a Dirichlet boundary condition is given on $\Gamma _{\text{out}} = \{762\} \times [0,15.24]$: $\Pres{l} = 10^6$ Pa, $\Sat{l} = 1$, and $\DensityComp{l}{h} = 0$. The upper and lower boundary is impervious, i.e. $\Flux{w} \cdot \nu = 0$ and $\Flux{h} \cdot \nu = 0$. Initial conditions are $\Pres{l} = 10^6$ Pa, $\Sat{l} = 1$, and $\DensityComp{l}{h} = 0$ for the whole domain. For the spatial discretization, we use a $100 \times 20$ mesh. The initial time step $dt = 20$ days and the end time is $T_{\text{final}} = 1160$ days. The initial smoothing parameter for the smooth Fischer-Burmeister function is $\tau = 10^{-6}$.
\par In the second example, the domain is a box of size $50m \times 30m \times 20m$. The porosity and permeability fields are random, generated by a geostatistic model using the open-source code MRST \cite{Lie11}. The porosity has a range of $[0.002,0.1]$ and the permeability varies from $1.377 \cdot 10^{-20}$ to $2.117 \cdot 10^{-15}$. Pure hydrogen is injected through the boundary at a corner: $\Flux{w} \cdot \nu = 0$ and $\Flux{h} \cdot \nu = 5.57 \times 10^{-2} kg/m^2/year$, and a Dirichlet boundary condition is chosen on the opposite corner: $\Pres{l} = 10^6$ Pa, $\Sat{l} = 1$, and $\DensityComp{l}{h} = 0$. The rest of the boundary is impervious, i.e. $\Flux{w} \cdot \nu = 0$ and $\Flux{h} \cdot \nu = 0$. Initial conditions are $\Pres{l} = 10^6$ Pa, $\Sat{l} = 1$, and $\DensityComp{l}{h} = 0$ for the whole domain. For the spatial discretization, we use a uniform $50 \times 30 \times 20$ mesh. The initial time step $dt = 200$ days and the end time is $T_{\text{final}} = 2000$ days. The initial smoothing parameter for the smooth Fischer-Burmeister function is $\tau = 10^{-4}$.
\begin{table}
\centering
{\renewcommand{\arraystretch}{1.2}
\begin{tabular}{l c c c c }
\cline{2-5}
\multicolumn{1}{c}{}& \multicolumn{2}{c}{Two-dimensional SPE10} & \multicolumn{2}{c}{Three-dimensional problem} \\
\hline
C-function &  FB & Smooth FB & FB & Smooth FB \\
\hline
Number of time steps & 60 (7) & 37 (4) & 13 (4) & 8 (3) \\
Average time step size (days) & 19.3 & 31.3 & 151.8 & 250.0\\ 
Total nonlinear iterations & 856 (147) & 530 (84) & 202 (84) & 135 (63)\\
Execution time (s) & 941.6 & 566.8 & 972.6 & 620.3 \\
\hline
\end{tabular}}
\captionsetup{justification=centering}
\caption{Performance comparison for heterogeneous problems.} \label{tab:result_heterogeneous}
\end{table}
\par For both of these problems, the semi-smooth Newton approach using the minimum function fails to converge for many time steps, and $dt$ becomes too small to perform the full simulation. Thus, only the results for the standard Fischer-Burmeister function and the smooth variant are reported in \cref{tab:result_heterogeneous}. Again, the numbers in parentheses are for the failed time steps and nonlinear iterations. The Jacobian smoothing method combined with the smooth Fischer-Burmeister function is more robust than the semi-smooth Newton approach with the standard Fischer-Burmeister function, as demonstrated by the reduction in the number successful and failed time steps. For example, in the two-dimensional SPE10 problem, the former requires only 37 successful time steps and registers 4 failed time steps, as opposed to 60 successful and 7 failed time steps of the latter. In terms of performance, the Jacobian smoothing method combined with the smooth Fischer-Burmeister function is clearly better as it helps decrease the number of nonlinear iterations and execution time by 34-40\% approximately for both the two-dimensional and three-dimensional problems.

\subsection{Scaling Results}\label{subsec:scaling}
To study parallel performance, we use the same setup as for the three-dimensional case with highly heterogeneous media considered in \cref{subsec:heterogeneous}. Parallel tests are run on Syrah, a Cray system with 5,184 Intel Xeon E5-2670 cores at the Lawrence Livermore National Laboratory Computing Center. Amanzi and other libraries are compiled with OpenMPI 1.6.5 and gcc-4.9.2. For strong scaling, the mesh size is fixed at $200 \times 120 \times 80$, and the problem has 5.76 million unknowns in total. We choose an initial time step of $dt = 2$ days and stop the simulation after 20 days. For weak scalability, the number of processors is increased in proportion to the problem size. We use meshes of size $50 \times 30 \times 20$, $100 \times 60 \times 40$, and $200 \times 120 \times 80$ with 2, 16, and 128 processors, respectively. The initial time step is set to $dt = 2$ days for all the mesh sizes and the simulation is stopped at $T = 200$ days. For both cases, the entry pressure is set at $\Pres{r} = 2\times 10^3$.
\begin{figure}
\centering
\includegraphics[width=0.75\textwidth]{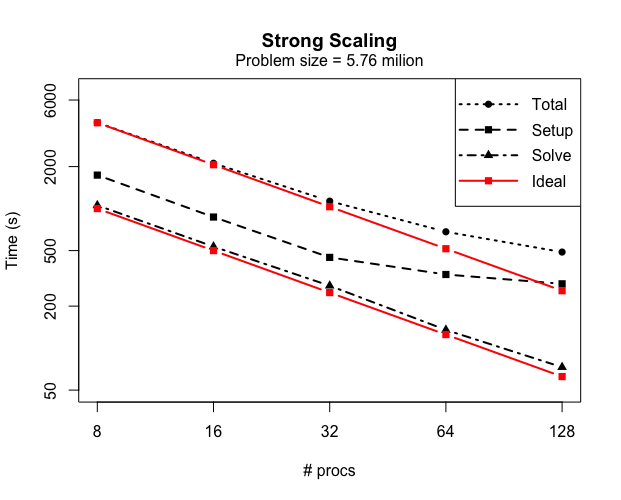}
\captionsetup{justification=centering}
\caption{Strong scaling for the three-dimensional heterogeneous problem. The total runtime for the simulation, the setup, and solve time for the linear solver are reported.}\label{fig:strong_scaling}
\end{figure}
The results reported in \cref{fig:strong_scaling} show that the Jacobian smoothing method, combined with GMRES preconditioned by hypreMGR achieves near optimal strong scalability on 8 to 128 processors for the total time needed to run the whole simulation. The slight deviation from the ideal performance at 64 and 128 processors results from the decrease in parallel performance of the setup phase of hypreMGR, which has been observed in \cite{Bui18}.
\par For weak scaling, a comparison between the Jacobian smoothing method using the smooth Fischer-Burmeister function and the semi-smooth Newton approach with the standard Fischer-Burmeister function is shown in \cref{tab:weak_scaling_jsm,tab:weak_scaling_fb}.
\begin{table}
\centering
{\renewcommand{\arraystretch}{1.2}
\begin{tabular}{l c c c }
\hline
Number of processors & 2 & 16 & 128 \\
Mesh size & $50 \times 30 \times 20$ & $100 \times 60 \times 40$ & $200 \times 120 \times 80$ \\
\hline
Initial smoothing parameter $\tau$ & $10^{-6}$ & $10^{-6}$ & $10^{-5}$ \\
Average step size (days) & 28.6 & 28.6 & 25.0\\
Number of time steps & 7 & 7 & 8\\
Average nonlinear iterations & 5.1 & 6.6 & 8.9\\
Average linear iterations & 10.7 & 13.5 & 17.4\\
Execution time & 122 (s) & 286 (s) & 995 (s) \\
\hline
\end{tabular}}
\captionsetup{justification=centering}
\caption{Weak scaling performance of the Jacobian smoothing method.} \label{tab:weak_scaling_jsm}
\end{table}
\begin{table}
\centering
{\renewcommand{\arraystretch}{1.2}
\begin{tabular}{l c c c }
\hline
Number of processors & 2 & 16 & 128 \\
Mesh size & $50 \times 30 \times 20$ & $100 \times 60 \times 40$ & $200 \times 120 \times 80$ \\
\hline
Average step size (days) & 28.6 & 25.0 & 3.45* \\
Number of time steps & 7 & 8 & 11 (2)* \\
Average nonlinear iterations & 4.7 & 6.6 & 11.1* \\
Average linear iterations & 12.7 & 22.0 & 28.5* \\
Execution time & 463 (s) & 1623 (s) & $>$ 4 hours \\
\hline
\end{tabular}}
\captionsetup{justification=centering}
\caption{Weak scaling performance of the semi-smooth Newton approach using the standard Fischer-Burmeister function.} \label{tab:weak_scaling_fb}
\end{table}
For the semi-smooth Newton method using the standard Fischer-Burmeister function, the simulation does not finish within the 4-hour limit of run time on the cluster. Thus, we only report the solver statistics up to $T = 38$ days when the simulation terminates. As the mesh is refined, the Jacobian smoothing method is clearly more robust and efficient than the semi-smooth Newton method using the standard Fischer-Burmeister function. Not only does it reduce the number of nonlinear iterations, it also helps improve the performance of the linear solver as indicated by smaller number of linear iterations. The execution time is significantly reduced as a consequence.

\section{Conclusions}\label{sec:conclusions}
In this work, we have developed a new Jacobian smoothing method based on the smooth Fischer-Burmeister function to solve the discrete nonlinear systems resulting from the the fully implicit discretization of the NCP formulation for compositional multiphase flow in porous media with phase transitions. Additionally, we performed various numerical experiments to compare our method with a semi-smooth Newton approach for two choices of C-function: the minimum and the Fischer-Burmeister functions. The results demonstrate that this method is significantly more robust and efficient with respect to the run time and number of nonlinear iterations. Unlike the semi-smooth Newton method using the minimum function, the Jacobian smoothing approach converges in all examples. Moreover, depending on the problem, it also reduces the number of nonlinear iterations and execution time by 34-40\% compared to the semi-smooth Newton method using the standard Fischer-Burmeister function.

\bibliographystyle{plainurl}
\bibliography{master_all}

\end{document}